\documentclass{amsart}
\usepackage{graphics}
\usepackage{ae}
\usepackage{amsmath,bbm}
\usepackage{dsfont}

\usepackage{amssymb}

\newcommand{\RR}{{\mathbb R}}

\newcommand{\NN}{{\mathbb N}}
\newcommand{\ZZ}{\mathbb Z}

\def\B{{\mathcal B}}

\def\P{{\mathcal P}}

\def\v{{\mathbf v}}

\def\x{{\mathbf x}}

\def\z{{\mathbf z}}
\def\1{\mathds{1}}
\def\Bbeta{J}

\numberwithin{equation}{section}

\newtheorem{theo}{Theorem}
\newtheorem*{theo*}{Theorem}
\newtheorem{prop}[theo]{Proposition}
\newtheorem{conj}[theo]{Conjecture}
\newtheorem{coro}[theo]{Corollary}
\newtheorem{lemma}[theo]{Lemma}
\theoremstyle{definition}
\newtheorem{defi}[theo]{Definition}
\newtheorem{rema}[theo]{Remark}
\newtheorem{exam}[theo]{Example}

\begin{document}

\title[Multidimensional extension of the Morse--Hedlund theorem]{Multidimensional extension of the Morse--Hedlund theorem}
\author{Fabien Durand}
\address[F.D.]{\newline
Universit\'e de Picardie Jules Verne\newline
Laboratoire Ami\'enois de Math\'ematiques Fondamentales et
Appliqu\'ees\newline
UMR 6140 CNRS UPJV\newline
33 rue Saint Leu\newline
80039 Amiens Cedex 01\newline
France.}
\email{fabien.durand@u-picardie.fr}
\author{Michel Rigo}
\address[M.R.]{\newline
Universit\'e de Li\`ege\newline
Institut de Math\'ematique\newline
Grande traverse 12 (B 37)\newline
B-4000 Li\`ege\newline
Belgium.}
\email{M.Rigo@ulg.ac.be}

\begin{abstract}
    A celebrated result of Morse and Hedlund, stated in 1938, asserts
    that a sequence $x$ over a finite alphabet is ultimately periodic
    if and only if, for some $n$, the number of different factors of
    length $n$ appearing in $x$ is less than $n+1$.  Attempts to
    extend this fundamental result, for example, to higher dimensions,
    have been considered during the last fifteen years.  Let $d\ge 2$.
    A legitimate extension to a multidimensional setting of the notion
    of periodicity is to consider sets of $\ZZ^d$ definable by a first
    order formula in the Presburger arithmetic
    $\langle\ZZ;<,+\rangle$.  With this latter notion and using a
    powerful criterion due to Muchnik, we exhibit a complete extension
    of the Morse--Hedlund theorem to an arbitrary dimension $d$ and
    characterize sets of $\ZZ^d$ definable in $\langle\ZZ;<,+\rangle$
    in terms of some functions counting recurrent blocks, that is,
    blocks occurring infinitely often.
\end{abstract}

\maketitle

\section{Introduction}
Let $d\ge 1$ be an integer. A set $M\subseteq \ZZ^d$ is said to be
{\em periodic} if there exists a non-zero vector
$\mathbf{p}\in\mathbb{Z}^d$ such that for all
$\mathbf{x}\in\mathbb{Z}^d$, $\mathbf{x}$ belongs to $M$ if and only
if $\mathbf{x}+\mathbf{p}$ belongs to $M$.  Periodic sets in dimension
greater than one have been investigated recently in the spirit of the
celebrated Morse--Hedlund theorem.

\begin{theo}[Morse--Hedlund theorem  \cite{Morse&Hedlund:1938}]
\label{th:morse-hedlund}
Let $A$ be a finite alphabet. Let $x\in A^{\mathbb N}$. Denote by
$p_x(n)$ the number of factors of length $n$ appearing in $x$.  Then
the following assertions are equivalent:
\begin{itemize}
  \item $x$ is ultimately periodic;
  \item there exists $n$ such that $p_x(n)\leq n$;
  \item there exist $n_0$ and $C$ such that $p_x(n)\leq C$ for all
    $n\geq n_0$.
\end{itemize}
\end{theo}

Observe that this theorem also can be transposed to any subset $M$ of
$\NN$ by considering its characteristic sequence over the alphabet $\{
0,1\}$.  The ultimate periodicity of this sequence is equivalent, for
$M$, to be a finite union of arithmetic progressions.

We are interested in generalizations to higher dimensions, that is, to
find conditions on the block counting function forcing (a notion
of) periodicity.  

Let us briefly mention some interesting works going in that direction.
We were first motivated by {\em Nivat's Conjecture}
\cite{Nivat:1997} stated in 1997, see also \cite{Berthe&Vuillon:2000},
asserting that if for some pair $(n_1,n_2)$ the number of 
$(n_1\times n_2)$-blocks appearing in $M\subset \ZZ^2$ is less or
equal to $n_1n_2$, then $M$ is periodic.  With the terminology
introduced in \cite{Sander&Tijdeman:2000} it would show that there is
a (restricted) {\em Periodicity Principle} in $\ZZ^2$.  This is also
strongly related to the {\em Period Conjecture} stated in
\cite{Lagarias&Pleasants:2002}.  More details on Nivat's conjecture,
Period Conjecture and related works can be found in
Section~\ref{section:comments}.

In this paper, we exhibit, thanks to Theorem~\ref{theo:main} given
below, a complete extension of the Morse--Hedlund theorem to an
arbitrary dimension $d$ and characterize sets of $\ZZ^d$ definable in
$\langle\ZZ;<,+\rangle$ in terms of some functions counting recurrent
blocks, that is, blocks occurring infinitely often. Note that we do
not consider ``purely'' periodic sets in $\ZZ^2$ or even $\ZZ^d$, with
$d\ge 2$, i.e., we are not looking at sets $M$ having a vector
$\mathbf{p}$ acting as ``global'' period. But instead we study what we
consider as a natural extension of the notion of periodicity to a
multidimensional setting.  Indeed, when A. L. Semenov extended in
\cite{Semenov:1977} the theorem of Cobham \cite{Cobham:1969} to a
multidimensional setting, it turns out that the notion of periodicity
he used was not the pure periodicity but sets in $\ZZ^d$ defined by a
first order formula of the Presburger arithmetic
$\langle\ZZ;<,+\rangle$.  It is important to observe that the subsets
of $\NN$ which are defined by a first order formula in
$\langle\NN;=,+\rangle$ are exactly the ultimately periodic sets
(i.e., finite union of arithmetic progressions). As an example, if
$$\varphi(x):=(x=3)\vee ((\exists y)(x=2y))\vee ((\exists
y)(x=5y+1))$$ is considered as a formula in $\langle\NN;=,+\rangle$,
then the set $$M=\{x\in\NN\mid \langle\NN;=,+\rangle\models
\varphi(x)\}=\{0,1,2,3,4,6,8,10,11,12,14,16,18,20,21,\ldots\}$$ is
ultimately periodic: for $x\ge 4$, $x\in M$ if and only if $x+10\in
M$. For a presentation of the sets in $\NN^d$ definable in the
Presburger arithmetic $\langle\NN;=,+\rangle$, see for instance the
very nice survey \cite{Bruyere&Hansel&Michaux&Villemaire:1994}. Except
for the scope of the variables, there is no fundamental difference
between the structures $\langle\NN;=,+\rangle$ and
$\langle\ZZ;<,+\rangle$. 
Therefore, in this paper, for the sake of simplicity our results are
described for the structure $\langle\ZZ;<,+\rangle$.  
Notice that in $\langle\NN;=,+\rangle$, the relation $x\le y$ can be
easily defined by $(\exists z)(x+z=y)$ but this is no more true when
$x,y,z$ are possibly negative and this is the reason why the order
relation $<$ has to be added to the Presburger arithmetic over $\ZZ$.
Note that over $\langle\ZZ;<,+\rangle$ the formula $x=y$ can be
defined by $\neg (x<y) \wedge \neg (y<x)$ and that $x=0$ can be
defined by $x+x=x$.

To obtain a characterization of sets $M$ in $\ZZ^d$ defined by a first
order formula of Presburger arithmetic, we estimate the function
$R_M(n)$ of different blocks of size $n$ occurring infinitely often in
$M$.  Collecting several results and observations about Presburger
arithmetic and mainly applying the powerful criterion of Muchnik
\cite{Muchnik:2003}, we derive the following theorem.
 
\begin{theo}\label{theo:main}
    Let $M$ be a subset in $\ZZ^d$. This set is definable in
    $\langle\ZZ;<,+\rangle$ if and only if
    $R_M(n)\in\mathcal{O}(n^{d-1})$ and every section is definable in
    $\langle\ZZ;<,+\rangle$.
\end{theo}

Hence, we use a different notion of periodicity as the one occurring
in Nivat's conjecture and a different complexity function but what we
obtain is a necessary and sufficient condition. 
Note that in many recent papers
\cite{Lagarias&Pleasants:2002,Nivat:1997,Quas&Zamboni:2004,Sander&Tijdeman:2000},
the given conditions are sufficient but are not necessary (there are
counter examples, see again Section~\ref{section:comments}).

If we apply Theorem~\ref{theo:main} recursively, we can express definability in
Presburger arithmetic solely in terms of recurrent block complexity.
More precisely, a set $M\subseteq\ZZ^d$ is definable in
$\langle\ZZ;<,+\rangle$ if and only if $R_M(n)\in\mathcal{O}(n^{d-1})$
and for all $k \in \{ 1, \dots , d-1 \}$ every $(d-k)$-dimensional
section has a recurrent block complexity in $\mathcal{O}(n^{d-k-1})$.

In Section~\ref{sec:1}, we define the terminology.
Section~\ref{sec:2} contains several examples of definable and
non-definable sets in $\ZZ^2$ to help the reader to develop some
intuition about Presburger arithmetic and complexity.  In Section
~\ref{sec:proofmain1} we prove the necessary condition relying mainly
on the elimination of quantifiers in Presburger arithmetic.  Observe
that for $d=1$ our main result corresponds to the Morse--Hedlund
theorem: a set $M\subseteq\mathbb{Z}$ is periodic, i.e., there exists
$t$ such that, for all $x\in\mathbb{Z}$, $x\in M\Leftrightarrow x+t\in
M$, if and only if there exists $n$ such that $p_M(n)\le n$
\cite{Morse&Hedlund:1938}.  For $d=2$ we give, in the first part of
Section~\ref{sec:proofmain2}, a direct proof using a lemma already
appearing in \cite{Epifanio&Koskas&Mignosi:2003} and
\cite{Quas&Zamboni:2004}.  In the second part of
Section~\ref{sec:proofmain2}, we treat the general multidimensional
case relying on completely different techniques developed by Muchnik
in \cite{Muchnik:2003}.  In Section \ref{section:comments} we describe
some works around the Nivat's Conjecture, the Periodicity Principle
and the Period Conjecture.  Moreover, under the extra hypothesis of
uniform recurrence, we obtain partial answers to these questions for
subsets of $\ZZ^d$.

\section{Setting up the framework}\label{sec:1}

Let $A$ be a non-empty finite alphabet and $d$ be a positive integer.  Let
$\mathbf{x}=(\mathbf{x}_1,\ldots,\mathbf{x}_d)\in\mathbb{Z}^d$. 
We consider
the norm $$||\mathbf{x}||=\sup_{1\le i\le d} |\mathbf{x}_i|.$$

  Let $k>0$. We define the {\em $k$-neighborhood} of size $k$ centered
  at $\mathbf{x}$ by
  \begin{equation}
      \label{eq:neighborhood}
      \B(\mathbf{x},k):=\{\mathbf{y}\in\ZZ^d\mid
  ||\mathbf{x}-\mathbf{y}||<k\}.
  \end{equation}
  Let $S\in A^{\ZZ^d}$ be an
  infinite $d$-dimensional word over $A$. It is convenient to consider
  $S$ as a map $S:\ZZ^d\to A$.  Let $S_{\mathbf{x},k}$ denote the {\it
    finite block of size $k$} given by the restriction of $S$ to the
  domain $\{\mathbf{x}\}+[\![0,k-1]\!]^d$. In this paper, when using
  notation like $[\![i,j]\!]$ with $i<j$, it has to be understood that
  we consider the set of integers $\{i,i+1,\ldots,j\}$. If $k=1$, we
  write $S_{\mathbf{x}}$ instead of $S_{\mathbf{x},1}$ to denote the
  letter in $S$ pointed at $\mathbf{x}$. To compare blocks, the
  respective domains of $S_{\mathbf{x},k}$ and $S_{\mathbf{y},k}$ can
  both be identified with $[\![0,k-1]\!]^d$ and the two maps can be
  compared.  Otherwise stated, we say that two blocks
  $S_{\mathbf{x},k}$ and $S_{\mathbf{y},k}$ of size $k$ are equal if
  and only if the mappings are the same, i.e., for all $\mathbf{v}$
  such that $\mathbf{v}\in[\![0,k-1]\!]^d$, we have
  $S_{\mathbf{x}+\mathbf{v}}=S_{\mathbf{y}+\mathbf{v}}$.

  The {\em block complexity} of $S$ is the map $p_S$ counting for each
  $n$ the number of different blocks $S_{\mathbf{x},n}$ appearing in
  $S$, i.e.,
  $$p_S:\NN\to\NN:n\mapsto\#\{S_{\mathbf{x},n}\mid
  \mathbf{x}\in\ZZ^d\}.$$
  The {\em recurrent block complexity} of $S$
  is the map $R_S$ counting for each $n$ the number of different
  blocks $S_{\mathbf{x},n}$ appearing infinitely often in $S$.
These blocks are called {\em recurrent blocks}.
Observe that we do not require to see such blocks in all quadrant.
In
  other words, this latter function counts only the different blocks
  $B$ of size $n$ such that for all $L\in\NN$, there exists
  $\mathbf{x}\in\ZZ^d$ verifying $||\mathbf{x}||\ge L$ and
  $S_{\mathbf{x},n}=B$. The following inequalities are obvious
  $$\forall n\ge 1,\quad R_S(n)\le p_S(n) \le (\# A)^{n^d}.$$

  \begin{rema}\label{rema:word-set}
      If the alphabet $A$ is binary, say $A=\{0,1\}$, then words in
      $A^{\ZZ^d}$ can be identified with subsets of $\ZZ^d$ and the
      definitions and notation given above can be applied to subsets
      of $\ZZ^d$. For instance, we can speak of the block complexity
      and recurrent block complexity of a subset $M$ of $\ZZ^d$. We
      denote these complexities respectively by $p_M$ and $R_M$.
  \end{rema}

\begin{rema}
Let $M = \{ 0 \}\subset \mathbb{Z}^d$. 
Then $R_M (n) = 1$ and $p_M (n) = n^d+1$ for all $n$.
With this example, we stress the fact that it is required to consider $R_M$ instead of $p_M$ in our main result.
\end{rema}

For the sake of completeness, we recall some essential facts about
sets definable by first order formulas in the Presburger arithmetic
$\langle\ZZ;<,+\rangle$. For more details, we refer the reader to
\cite{Bes:2001,Bruyere&Hansel&Michaux&Villemaire:1994,Muchnik:2003}.
\begin{rema}\label{rema:d=1}
Recall that the subsets of $\NN$ definable in
$\langle\NN;=,+\rangle$ are exactly the ultimately periodic sets in
$\NN$. When taking a subset $M$ of $\ZZ$ definable in
$\langle\ZZ;<,+\rangle$, we can consider the two subsets of $\NN$,
$$M\cap\ZZ_{\ge 0} \text{ and }-(M\cap\ZZ_{<0})=\{-x\mid (x\in M)
\wedge (x<0)\}$$
which are both definable in $\langle\NN;=,+\rangle$
and therefore ultimately periodic. Consequently, if $M\subseteq\ZZ$ is
definable in $\langle\ZZ;<,+\rangle$, there exist $N,p,q>0$ (which are
chosen minimal) such that for $x>N$, $x\in M$ if and only if $x+p\in
M$ and for $x<-N$, $x\in M$ if and only if $x-q\in M$. Notice that
even in dimension one, if $p\neq q$ or $N\neq 0$, then $M$ which is definable in
$\langle\ZZ;<,+\rangle$ cannot be ``purely'' periodic. 
As a trivial example, consider the set $M=\{0\}\subset\mathbb{Z}$.
Then, $p_M(n)=n+1$ for all $n\ge 1$.
\end{rema}

In this context, Presburger definability can be viewed as a natural
extension of ultimate periodicity in dimensions higher than one.
Before stating the so-called Muchnik's criterion, we need a few definitions.
\begin{defi}
    A set $M\subseteq\ZZ^d$ is said to be {\em linear}, if there exist
    $\mathbf{x}\in\ZZ^d$ and a finite set (possibly empty)
    $V\subset\ZZ^d$ such that $$M=\{\mathbf{x}\}+\sum_{\mathbf{v}\in
      V}\NN \mathbf{v}.$$ A finite union of linear sets is a {\it
      semi-linear} set.
\end{defi}

Let $M\subseteq\ZZ^d$, $i\in\{1,\ldots,d\}$ and $c\in\ZZ$. The {\it
  section} $M_{i,c}\subseteq\ZZ^{d-1}$ is given by
$$M_{i,c}=\{(x_1,\ldots,x_{i-1},x_{i+1},\ldots,x_d)\mid
(x_1,\ldots,x_{i-1},c,x_{i+1},\ldots,x_d)\in M\}.$$

\begin{defi}
    Let $\mathbf{v}\in\ZZ^d\setminus\{\mathbf{0}\}$ and
    $M,X\subseteq\ZZ^d$. The set $M$ is said to be {\it
      $\mathbf{v}$-periodic inside} $X$ if for any
    $\mathbf{m},\mathbf{m}+\mathbf{v}\in X$,
    $$\mathbf{m}\in M\Leftrightarrow \mathbf{m}+\mathbf{v}\in M$$
    (i.e., any two points inside X that differ by $\mathbf{v}$ either
    both belong to $M$ or both do not belong to $M$).
We say $\mathbf{v}$ is a {\em local period} for $X$.
    When $V$ is a subset of $\mathbb{Z}^d$, we say $M$ is {\em $V$-periodic inside} $X$ if $M$ is $\mathbf{v}$-periodic inside $X$ for some $\mathbf{v}\in V$.
    \begin{figure}[htbp]
      \centering
      \scalebox{.75}{\includegraphics{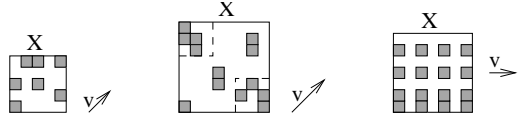}}
      \caption{Three examples of $\mathbf{v}$-periodicity inside $X$.}
      \label{fig:10}
  \end{figure}
\end{defi}

As we will see definability in $\langle \mathbb{Z} ; < , + \rangle $ can be expressed in terms of 
local periodicity as given in \cite{Bruyere&Hansel&Michaux&Villemaire:1994}.

\begin{defi}\label{def:localp}
    A set $M\subseteq\ZZ^d$ is said to be {\em locally periodic} if
    there exists a finite set $V\subset\ZZ^d\setminus\{\mathbf{0}\}$
    such that for some $K>\sum_{\mathbf{v}\in V} ||\mathbf{v}||$ and
    $L\ge 0$, $$\forall \mathbf{x}\in\ZZ^d, ||\mathbf{x}||\ge L,
    \exists \mathbf{v}\in V: M \text{ is }\mathbf{v}\text{-periodic
      inside }\B(\mathbf{x},K).$$    
    Figure~\ref{fig:new01} illustrates this notion: far enough of the
    origin, the set $M$ is $\mathbf{v}$-periodic inside any
    $K$-neigborhood for some $\mathbf{v}\in V$.
\begin{figure}[htbp]
      \centering
      \scalebox{.75}{\includegraphics{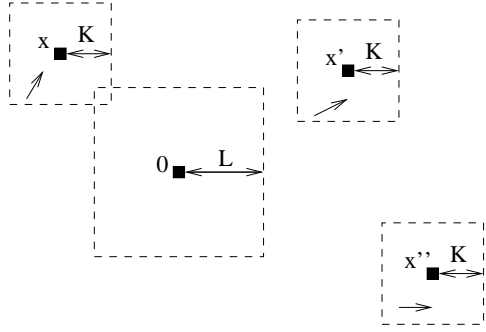}}
      \caption{Illustration of local periodicity.}
      \label{fig:new01}
  \end{figure}
\end{defi}

In \cite{Muchnik:2003} the following alternative definition of local periodicity is given,
see also \cite{Bes:2001,Michaux&Villemaire:1996}. 
Moreover, in \cite{Quas&Zamboni:2004} the authors also defined a related notion of local periodicity.

\begin{defi}
    Let $M\subseteq\ZZ^d$.  
    If there exists  a finite subset $V$ of $\ZZ^d \setminus\{\mathbf{0}\}$ verifying: for all
  $K\ge 1$ there exists $L\in\NN$ such that
    $$
    \forall \mathbf{x}\in\ZZ^d, ||\mathbf{x}||\ge L,
    \exists \mathbf{v}\in V: M \text{ is }\mathbf{v}\text{-periodic
      inside }\B(\mathbf{x},K),
    $$
    then we say that $M$ {\em satisfies Muchnik's
    condition}. 
A picture similar to the one in  Figure~\ref{fig:new01} should also illustrate this notion.
\end{defi}

As shown by the following theorem, these two definitions play similar role to characterize definability. The reader may already notice that the sections of $M$ enter the picture.

\begin{theo}\cite{Muchnik:2003}
\label{theo:muchnik}
    Let $M$ be a subset of $\ZZ^d$. The following statements are
    equivalent.
    \begin{enumerate}
      \item $M$ is definable in $\langle\ZZ;<,+\rangle$,
      \item $M$ is semi-linear,
      \item \label{theo:equivlocper} every section of $M$ is definable in
        $\langle\ZZ;<,+\rangle$ and $M$ is locally periodic,
      \item \label{theo:muchnikcond} every section of $M$ is definable in
        $\langle\ZZ;<,+\rangle$ and $M$ satisfies Muchnik's condition.
    \end{enumerate}
\end{theo}

\begin{rema}

    It is clear in Theorem~\ref{theo:muchnik} that condition \eqref{theo:muchnikcond}
    implies condition \eqref{theo:equivlocper}, i.e., Muchnik's condition implies local
    periodicity. The converse is not so obvious and as observed in
    \cite[p. 1438]{Muchnik:2003}, a close inspection of the proof of Theorem~\ref{theo:muchnik} reveals that assuming \eqref{theo:equivlocper} is enough to derive
the
    definability of $M$ and that this definability implies Muchnik's
    condition.
\end{rema}

Let $A$ be a finite alphabet, $d\ge 1$ and $S\in A^{\ZZ^d}$ be a
    $d$-dimensional word over $A$.  For each $a\in A$, we define the set
    $$S^{-1}(a):=\{\mathbf{x}\in\ZZ^d \mid S_{\mathbf{x}}=a\}.$$
    
\begin{defi}
    The
    word $S$ is {\em definable} in $\langle\ZZ;<,+\rangle$ if and only
    if for all $a\in A$, the sets $S^{-1}(a)$ are definable in
    $\langle\ZZ;<,+\rangle$. 
\end{defi}

We can now recall our main result (stated in the introduction as Theorem~\ref{theo:main}).
    \begin{theo*}
        Let $M$ be a subset in $\ZZ^d$. This set is definable in
        $\langle\ZZ;<,+\rangle$ if and only if
        $R_M(n)\in\mathcal{O}(n^{d-1})$ and every section is definable
        in $\langle\ZZ;<,+\rangle$.
    \end{theo*}

    The case $d=1$ is equivalent to Morse--Hedlund theorem and was already discussed in
    Remark~\ref{rema:d=1}. Notice that the number of recurrent
    factors occurring in an ultimately periodic word in $\{0,1\}^\NN$
    is bounded by a constant. The proof of Theorem~\ref{theo:main} is given in
    Section~\ref{sec:proofmain1} and Section~\ref{sec:proofmain2}. But prior to these
    developments, let us make a few immediate comments.
    Theorem~\ref{theo:main} has an immediate application.

    \begin{coro}
    \label{coro:1}
        Let $A$ be a finite alphabet and $d\ge 1$. A $d$-dimensional
        word $S\in A^{\ZZ^d}$ is definable in $\langle\ZZ;<,+\rangle$
        if and only if $R_S(n)\in\mathcal{O}(n^{d-1})$ and every section
        is definable in $\langle\ZZ;<,+\rangle$.
    \end{coro}

Let $k \in \{ 1, \dots ,d-1 \}$. 
A {\em $(d-k)$-dimensional section} of $M\subseteq\ZZ^d$ is
a subset of $\ZZ^{d-k}$ where $k$ components have been fixed to
constants. So we define the set
$$M_{(i_1,\ldots,i_k),(c_1,\ldots,c_k)}=\{ \mathbf{x}\in M \mid
\mathbf{x}_j=c_j,\quad j=1,\ldots,k\}$$
viewed as a subset of
$\ZZ^{d-k}$. 
If we apply Theorem~\ref{theo:main} recursively, we get
the following result expressing definability in Presburger arithmetic
solely in terms of recurrent block complexity.

\begin{coro}
\label{coro:2}
    A set $M\subseteq\ZZ^d$ is definable in $\langle\ZZ;<,+\rangle$
    if and only if $R_M(n)\in\mathcal{O}(n^{d-1})$ and, for all $k\in\{1,\ldots,d-1\}$, every $(d-k)$-dimensional section has a recurrent block
    complexity in $\mathcal{O}(n^{d-k-1})$.
\end{coro}

Using similar arguments as in Remark~\ref{rema:d=1} we can replace in Theorem~\ref{theo:main}, Corollary~\ref{coro:1} and Corollary~\ref{coro:2}, $\ZZ^d$ by $\NN^d$ and $\langle\ZZ;<,+\rangle$ by $\langle\NN;=,+\rangle$.

\section{Running examples}\label{sec:2}

We briefly give some examples that we hope could be helpful.

\begin{exam}\label{exa:1}
    Consider the set $M\subset\ZZ^2$ represented on the left in
    Figure~\ref{fig:03} and defined by the formula
  $$\varphi(x,y):= (x\ge 0)\wedge (y\ge 0)\wedge (\exists
\lambda)((x,y)=\lambda
  (1,1))\wedge(\exists \lambda)((x,y)=(0,1)+\lambda (1,0)).$$

  \begin{figure}[htbp]
    \centering
    \scalebox{0.75}{\includegraphics{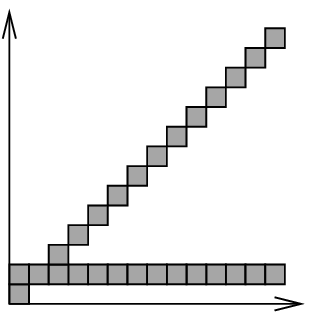}}\quad\quad
\scalebox{0.75}{\includegraphics{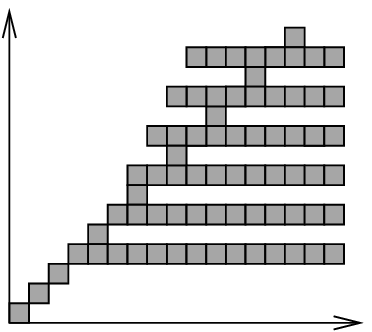}}
    \caption{Examples of sets definable in Presburger arithmetic.}
    \label{fig:03}
  \end{figure}

  One can show that the block complexity $p_M$ satisfies $p_M(n)\ge n^2$ for $n$
  large enough. A simple argument is to observe that the $n^2$ blocks
  of size $n$, $(i,j)+[\![0,n-1]\!]^2$ for $i,j\in[\![-n+1,0]\!]$ are
  all distinct.  Local periodicity is clear with a set of periods
  equal to $V=\{(1,1),(1,0)\}$, and, $K=3$ and $L=4$. Muchnik's
  condition is also obvious to check with the same set of periods. The
  recurrent block complexity is $R_M(n)=3\, n$. 
\end{exam}

\begin{exam}\label{exa:2}
    Consider the set $M\subset\ZZ^2$ represented on the right in
    Figure~\ref{fig:03} and defined by the formula
    \begin{eqnarray*}
    \varphi(x,y)&:=& (x\ge 0)\wedge (y\ge 0)\wedge ((\exists
\lambda)((x,y)=\lambda
    (1,1))\\
&&\vee(\exists \lambda)(\exists\mu)((x,y)=(4,3)+\lambda
    (1,0)+\mu (1,2))).
    \end{eqnarray*}
    Notice that this set can also be defined with
    a quantifier-free formula
    $\psi(x,y)=(x\ge 0)\wedge (y\ge
0)\wedge(\psi_1(x,y)\vee\psi_2(x,y)\vee\psi_3(x,y))$ where
\begin{equation}
    \label{eq:quantif}
    \begin{array}{l}
\psi_1(x,y):= (y\ge x) \wedge (2x\ge y+5) \wedge (y\equiv 1\mod{2}),\\
\psi_2(x,y):= (y\ge x)\wedge (x\ge y),\\
\psi_3(x,y):= (x\ge y)\wedge (y\ge 3)\wedge (y\equiv 1\mod{2}).
    \end{array}
\end{equation}
    It is not difficult to compute that the
    recurrent block complexity is $R_M(n)=7\, n-1$ for all $n\ge 1$.
    \end{exam}

The following two examples will be used to see that our main result
cannot be improved.

\begin{exam}[Fibonacci]
    Let $\mathcal{F}=(f_i)_{i\ge 0}=0100101001\cdots$ be the Fibonacci
    word generated by the morphism $h:0\mapsto 01, 1\mapsto 0$. It is
    well-known that $p_{\mathcal{F}}(n)=n+1$. Now consider the word
$\mathcal{G}=(g_i)_{i\in\ZZ}$ indexed by $\ZZ$ and defined by
$$g_i=\left\{\begin{array}{ll}
g_i=f_i&\text{ if }i\ge 0,\\
g_i=1&\text{ otherwise}\\
\end{array}\right.$$
$$g=\cdots 11111111.0100101001\cdots .$$
We have
$p_{\mathcal{G}}(n)=2n$.  Consider the set $G\subset{\ZZ^d}$ defined
by $\mathbf{x}\in G$ if and only if 
$\mathbf{x} = (\mathbf{x}_1, \dots , \mathbf{x}_d)$ is such that $g_{\mathbf{x}_1}=1$.  
It is clear that $G$ is locally periodic and that 
$R_G(n)=p_G(n)=2n$.
None of the sections $G_{j,c}$ for $j\ge
2$ can be defined in $\langle\ZZ;<,+\rangle$. Indeed, the Fibonacci
word is not ultimately periodic (see \cite{Morse&Hedlund:1940}) and therefore not definable in
$\langle\NN;=,+\rangle$.

Thus, the hypothesis on the sections cannot be remove.
It also shows that it is not enough to consider $p_G$ instead of $R_G$.
\end{exam}

\begin{exam}[Toeplitz'like set]
  Consider the set $T\subset\ZZ^2$ depicted in Figure~\ref{fig:02}
  such that
  $$(i,j)\in T \Leftrightarrow (i,j\ge 0) \wedge (\exists n>0: i=n\,
2^{j+1}).$$
  \begin{figure}[htbp]
    \centering
    \scalebox{0.6}{\includegraphics{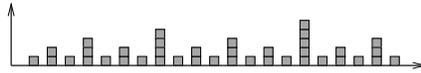}}
    \caption{Toeplitz'like set.}
    \label{fig:02}
  \end{figure}
  Each section is definable in $\langle\ZZ;<,+\rangle$ but the set is
  not locally periodic and does not satisfy Muchnik's condition. 
  It is
  not difficult to see that $R_T(n)\ge n^2$ since the patterns
  represented in Figure~\ref{fig:05} appear infinitely often.
  Thus, the local periodicity property cannot be remove.
  \begin{figure}[htbp]
      \centering
      \scalebox{0.6}{\includegraphics{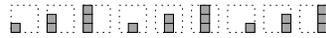}}
      \caption{Some blocks of size $3$ occurring in $T$.}
      \label{fig:05}
  \end{figure}
  We call it {\em Toeplitz'like set} because, one can obtain the
  heights of the successive ``heaps'' as follows. At first step,
  consider the infinite word $(0?)^\omega=0?0?\cdots$, on the second
  step, replace in order all the ``$?$'' symbols with the letters of
  the infinite word $(1?)^\omega$. On the $(n+1)$-th step, replace the
  remaining ``$?$'' with the letters of the word $(n?)^\omega$. This
  process generates an infinite word over the alphabet $\NN$ whose
  $j$-th element is the height of the $j$-th heap.
\end{exam}

\section{Complexity of Presburger definable sets}
\label{sec:proofmain1}

In what follows, as explained in Remark~\ref{rema:word-set}, we make
no distinction between a subset of $\ZZ^d$ and its characteristic word
which belongs to $\{0,1\}^{\ZZ^d}$. We consider indifferently
complexity functions of both subsets of $\ZZ^d$ and binary words indexed by
$\ZZ^d$.
To count the number of different blocks of size $n$ it is convenient to use the following notation:

$$
\mathcal{C} (\mathbf{x} , n) = \{ \mathbf{y} \in \mathbb{Z}^d \mid \mathbf{y}=\mathbf{x}+\mathbf{z} , \mathbf{z}\in [\![0,n-1]\!]^d  \}. 
$$
One should avoid the confusion between the neighborhood $\mathcal{B}(\mathbf{x} , n)$ centered at $\mathbf{x}$ given in~\eqref{eq:neighborhood} and this set $\mathcal{C} (\mathbf{x} , n)$ {\em pointed} at $\mathbf{x}$.
\begin{prop}
  If $X\subseteq\ZZ^d$ is definable in $\langle\ZZ ; =,+\rangle$, then
  $R_X(n)\in\mathcal{O}(n^{d-1})$.
\end{prop}

For a better understanding of the main arguments, we first give the proof
for $d=2$.

\begin{proof}
    Let $X\subseteq\ZZ^2$ be definable in $\langle\ZZ;=,+\rangle$. We
    can assume \cite{Presburger:1929,Presburger:1991} that the formula defining $X$ is a
    finite boolean combination of formulas of the kind
    $$\lambda_i\, x+\mu_i\, y\ge \nu_i ,\quad i\in\{1,\ldots,r\} \quad
    \text{ or }\quad \lambda_i\, x+\mu_i\, y \equiv \alpha_i
    \pmod{\Bbeta},\quad i\in\{1,\ldots,s\}$$
    where the
    $\alpha_i,\Bbeta,\lambda_i,\mu_i,\nu_i$'s are integers (for an
    example, see for instance formula \eqref{eq:quantif}). Therefore,
    except for a neighborhood of the origin, $X$ is made up of finite
    number of regions bounded by two half (straight) lines. Inside such
    a region $X$ is periodic (the period being determined with the
    constant $\Bbeta$ appearing in the formulas given above).
In figure~\ref{fig:01}, no two lines have the same slope. Let us make
    this assumption first.
  \begin{figure}[htbp]
    \centering
    \scalebox{.6}{\includegraphics{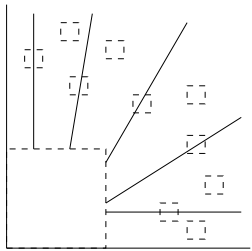}}
    \caption{Regions in dimension $2$.}
    \label{fig:01}
  \end{figure}

  Let $n>0$.  
Note that bounded polygonal regions delimited by several lines
  remain in a finite domain close to the origin. Hence, only infinite
  regions have to be taken into account to estimate $R_X$.  Thus we consider
  one infinite region bounded by two lines $\lambda_i\, x+\mu_i\, y=
  \nu_i$, $i=1,2$.  Since their slopes are different, there exists
  $s(n)$ such that for any $\mathbf{x}$ of norm $||\mathbf{x}||>s(n)$,
  $\mathcal{C}(\mathbf{x},n)$ intersects at most one line. The number
  of distinct blocks of size $n$ lying completely in one region $R$ is
  bounded by a constant $r_R$ (due to the periodicity given by a
  finite number of congruence relations). The number of different
  blocks of size $n$ intersecting two regions is bounded by $C\, n$
  where the constant $C$ depends on the periods within the two regions
  and the period on the line. Indeed the number of relative positions
  of a square of size $n$ with respect to a given line is bounded by
  $2n$ (see Figure~\ref{fig:13}).
  \begin{figure}[htbp]
    \centering
    \scalebox{1}{\includegraphics{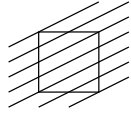}}
    \caption{Intersection of a square and a line.}
    \label{fig:13}
  \end{figure}
  For each such position, the number of different patterns within the
  two regions and on the line is bounded respectively by three
  constants $r_1,r_2,r_3$ (one for each region and one for each line). This gives a number of different blocks
  bounded by $2r_1r_2r_3n$.

  The conclusion follows from the fact that the number of regions is
  finite. 

  We have first considered the case where the (infinite) regions were
  defined by lines having distinct slopes.  Now, if several lines
  defining regions have the same slope as depicted in
  Figure~\ref{fig:new02}, then one can proceed using the same
  arguments.  \begin{figure}[htbp] \centering
      \scalebox{.6}{\includegraphics{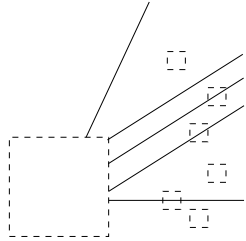}}
    \caption{Some regions defined by parallel lines.}
    \label{fig:new02}
  \end{figure}
  Parallel lines do not increase the linear bound on the complexity.
  Indeed, the number of relative positions of a square of size $n$
  with respect to given parallel lines is bounded by $D\, n$, for some
  $D$.

\medskip

Now consider the general case.

\medskip

    Thanks to quantifier elimination in Presburger arithmetic
    \cite{Enderton:2001}, we can assume that the formula defining $X$ is a
    finite disjunction of formulas $\varphi_i(x_1,\ldots,x_d)$,
    $i=1,\ldots,k$, where each such formula is a conjunction of the
    form

\begin{equation}
\label{canon-formula}
\varphi_i(x_1,\ldots,x_d):= \bigwedge_{j=1}^{\ell} t_{i,j}(x_1,\ldots,x_d)
\end{equation}

with
\begin{equation}
\label{domain}
t_{i,j}(x_1,\ldots,x_d):= \sum_{k=1}^d u^{(i,j)}_{k}\, x_k\ge
c_{i,j}
\end{equation}

or 

\begin{equation}
\label{pattern}
t_{i,j}(x_1,\ldots,x_d):=
\sum_{k=1}^d u^{(i,j)}_{k} x_k \equiv e_{i,j} \pmod{\Bbeta}
\end{equation}

where all the
occurring constants are integers.  Notice that we can choose the same
constant $\Bbeta$ in every modular equation occurring in the definition
of $X$.  
Each $\varphi_i$ defines a domain\footnote{We call {\em domain} any subset of $\ZZ^d$ defined by a conjunction of formulas of the kind \eqref{domain}.} of $\ZZ^d$
intersected with a pattern\footnote{We call {\em pattern} any subset of $\ZZ^d$ defined by a conjunction of formulas of the kind \eqref{pattern}.}. The following lemma is obvious.

\begin{lemma}
\label{lem:nonintersecting}
    The number of distinct blocks of size $n$ for a
    set $X\subset \ZZ^d$ defined by a boolean combination of formulas of
the kind
$\sum_{k=1}^d u_{k}\,
x_k \equiv e \pmod{\Bbeta}$ is bounded from above by $\Bbeta^d$.
\end{lemma}

Let $\pi_{i,j}$ be the hyperplane having equation $\sum
u^{(i,j)}_{k}\, x_k= c_{i,j}$ corresponding to~\eqref{domain}. We define $\pi_{i,j}'$ as the vector subspace
having equation $\sum u^{(i,j)}_{k}\, x_k= 0$. For the sake of
simplicity, let us enumerate all the distinct hyperplanes occurring in
the quantifier-free formula as $\Pi=\{\pi_1,\ldots,\pi_t\}$ (with the
corresponding vector subspaces denoted by $\pi_1',\ldots,\pi_t'$).

\medskip

Let $B$ be a recurrent block of size $n$, i.e., there exist
infinitely many $\mathbf{y}_1,\mathbf{y}_2,\ldots$ in $\mathbb{Z}^d$ such that
$$B=X_{\mathbf{y}_1,n}=X_{\mathbf{y}_2,n}=\cdots.$$
Our aim is to
estimate the number of distinct such blocks. 
Assume that for large enough $i$, $\mathcal{C}(\mathbf{y}_i,n)$ intersect\footnote{In what follows, when dealing with intersections, one should understand that the underlying hypercube, considered as a convex body in $\mathbb{R}^d$, does not intersect any hyperplane.}
no hyperplane in $\Pi$. Then it lies completely in a single domain. We
can apply the previous lemma. There are at most $\Bbeta^d$ different
such blocks for each domain and the number of domains is finite (and determined
by the formula).
\medskip

In what follows, we may therefore assume that all the blocks
 $X_{\mathbf{y}_i,n}$ intersect each of the
hyperplanes $\pi_1,\ldots,\pi_j$ and no other, the $\mathbf{y}_i$
belonging to the same domain $D$.
We can do the reasoning for any
subset of $\Pi$ and any domain $D$ (there are at most $ 2^t$ subsets of $\Pi$).
For the sake of simplicity, we have chosen the first $j$ elements in
$\Pi$.
Then, it suffices to prove that there exists a constant $C$, not depending on
$n$, such that the number
of pairwise distinct recurrent blocks $X_{\mathbf{y},n}$ intersecting each of the
hyperplanes $\pi_1,\ldots,\pi_j$ and no other, with $\mathbf{y}\in D$,
is bounded by $Cn^{d-1}$.

\begin{lemma}\label{lem:dimW}
    Let $B$ be a block of size $n$. If there exist infinitely many
    (different) vectors $\mathbf{y}_1,\mathbf{y}_2,\ldots\in\mathbb{Z}^d$ such that, 
    for all $i$, 
    $\mathcal{C}(\mathbf{y}_i,n)$ intersects each of the hyperplanes
    $\pi_1,\ldots,\pi_j$, then
$$\dim (\pi_1'\cap\cdots \cap\pi_j')\ge 1.$$
\end{lemma}

\begin{proof}
    Assume to the contrary that the intersection of the vector
    subspaces $\pi_1'\cap\cdots \cap\pi_j'$ is reduced to
    $\{\mathbf{0}\}$. Therefore, for a given $n$, there are finitely many
    hypercubes $\mathbf{z}+[0,n]^d$ of size $n$ intersecting
    simultaneously $\pi_1,\ldots,\pi_j$, with
    $\mathbf{z}\in\mathbb{Z}^d$.
\end{proof}

Now let us assume that $W=\pi_1'\cap\cdots \cap\pi_j'$ is a vector
space of dimension $\ell \in \{ 1 , \dots , d-1 \}$.
Let $(\mathbf{b}_1,\ldots,\mathbf{b}_\ell)$ be a basis of $W$.
We can suppose the entries of $\mathbf{b}_1,\ldots,\mathbf{b}_\ell$ are
integers and multiples of
$\Bbeta$.

The first idea is to show that a move parallel to $W$ only produces a constant number of distinct blocks.

\begin{lemma}\label{lem:11}
    Suppose $\mathcal{C}(\mathbf{v},n)$, $\mathbf{v} \in D$,  intersects each of the
hyperplanes
    $\pi_1,\ldots,\pi_j$ and no other.
Then for all $\mathbf{w}=\alpha_1\,
    \mathbf{b}_1+\cdots+\alpha_\ell\, \mathbf{b}_\ell$ with
    $\alpha_1,\ldots,\alpha_\ell\in\mathbb{Z}$,
    $\mathcal{C}(\mathbf{v}+\mathbf{w},n)$ intersects each of the hyperplanes
    $\pi_1,\ldots,\pi_j$.
    Furthermore, if $\mathbf{v}+\mathbf{w} \in D$ and $\mathcal{C}(\mathbf{v}+\mathbf{w},n)$ intersects no other hyperplane,
    then

$$
X_{\mathbf{v},n}=X_{\mathbf{v}+\mathbf{w},n}.
$$
\end{lemma}

\begin{proof}
    By definition of $W$, the hypercubes $\mathcal{C}(\mathbf{v},n)$
    and $\mathcal{C}(\mathbf{v}+\mathbf{w},n)$ have the same relative
    position with respect to $\pi_1,\ldots,\pi_j$ 
    (e.g., consider two intersecting planes in a space of dimension
    $3$ and assume that a cube intersects these two planes, then a
    translation of $\mathbf{w}$ corresponds to a move in the direction
    of the intersection of the two planes, so the cube still
    intersects these two planes with a similar relative position of
    the cube compared with the planes).  Moreover, since in each
    domain the pattern is given by congruence relations modulo
    $\Bbeta$, we conclude thanks to the choice made for the basis
    (each of the components is a multiple of $\Bbeta$).
\end{proof}

\begin{rema}\label{rem:11}
    There exist finitely many vectors of the kind $\alpha_1\,
    \mathbf{b}_1+\cdots+\alpha_\ell\, \mathbf{b}_\ell$ belonging
    $\ZZ^d$ with $(\alpha_1,\ldots,\alpha_\ell)\in\mathbb{Q}^\ell\cap
    [-1,1]^\ell$. This observation and the periodicity induced by the
    previous lemma lead to the following fact. Let $\mathbf{v} \in D$
    be such that $\mathcal{C}(\mathbf{v},n)$ intersects each of the
    hyperplanes $\pi_1,\ldots,\pi_j$ and no other. There exists a
    constant $Q$ (independent of $n$) such that the number of pairwise
    distinct blocks $X_{\mathbf{v}+\mathbf{w},n}$ is bounded by $Q$
    when $\alpha_1,\ldots,\alpha_\ell\in\mathbb{Q}$ are such that
    $\mathbf{w}=\alpha_1\, \mathbf{b}_1+\cdots+\alpha_\ell\,
    \mathbf{b}_\ell\in\ZZ^d$, $\mathbf{v}+\mathbf{w} \in D$ and
    $\mathcal{C}(\mathbf{v}+\mathbf{w},n)$ intersects each of the
    hyperplanes $\pi_1,\ldots,\pi_j$ and no other.
\end{rema}

Let $\mathbf{e}_1,\ldots,\mathbf{e}_{d-\ell}$ be integer vectors such
that
$(\mathbf{e}_1,\ldots,\mathbf{e}_{d-\ell},\mathbf{b}_1,\ldots,\mathbf{b}_\ell)$
is a basis of the space $\RR^d$. The idea is now to move in a direction not parallel to $W$ and show that we have $d-\ell$ degrees of freedom.

\begin{lemma}\label{lem:12}
    For all $i\in\{1,\ldots,d-\ell\}$, 
    there exists a constant $C_i>0$
    such that for all 
    $\mathcal{C}(\mathbf{x},n)$ intersecting $\pi_1,\ldots,\pi_j$ at some point $\mathbf{y}\in\RR^d$, the following holds:
    \begin{itemize}
      \item if $\mathcal{C}(\mathbf{x}+k\mathbf{e}_i,n)$ intersects
    $\pi_1,\ldots,\pi_j$, then $|k|\le C_i n$,
  \item if $|k|>C_i n$, then $\mathcal{C}(\mathbf{x}+k\mathbf{e}_i,n)$
    does not intersect $\pi_1,\ldots,\pi_j$.
    \end{itemize}
\end{lemma}

\begin{proof}
    Without loss of generality, we may assume $i=1$. Since $\mathbf{e}_1$ does not belong to $W$,
    there exists $m\in\{1,\ldots,j\}$ such that $\mathbf{e}_1$ does not belong to the vector subspace
    $\pi_m'$. In a system of coordinates built on the basis $(\mathbf{e}_1,\ldots,\mathbf{e}_{d-\ell},\mathbf{b}_1,\ldots,\mathbf{b}_\ell)$, we may assume $\pi_m$ has an equation of the kind 
    $$\alpha_1^{(m)} \, x_1+\cdots+\alpha_d^{(m)} \, x_d=c$$
    and that
    $\mathbf{y}=(a_1,a_2,\ldots,a_d)$ belongs to $\pi_m$. In this case, the equation of the corresponding vector subspace
    $\pi_m'$ is $\alpha_1^{(m)} \, x_1+\cdots+\alpha_d^{(m)} \, x_d=0$. The components of $\mathbf{e}_1$ in the basis $(\mathbf{e}_1,\ldots,\mathbf{e}_{d-\ell},\mathbf{b}_1,\ldots,\mathbf{b}_\ell)$ are $(1,0,\ldots,0)$ and they do not satisfy the equation of $\pi_m'$ (otherwise, we would get that $\mathbf{e}_1\in\pi_m'$). Therefore, $\alpha_1^{(m)}$ is non-zero.

We can determine
    the values of $k$ such that $(a_1+k,a_2+\ell_2,\ldots,a_d+\ell_d)$
    belongs to $\pi_m$ for some $\ell_2,\ldots,\ell_d\in [\![-n,n]\!]$.
    Indeed, this
    latter point belongs to $\pi_m$ if and only if
    $$\exists \ell_2,\ldots,\ell_d\in [\![-n,n]\!] : \alpha_1^{(m)}\, k+\alpha_2^{(m)}\,
    \ell_2+\cdots+\alpha_d^{(m)}\, \ell_d=0.$$
    That is, if and only if
    $$k=-(\alpha_2^{(m)}\, \ell_2+\cdots+\alpha_d^{(m)}\, \ell_d)/\alpha_1^{(m)}.$$
 
 Thus, it suffices to take 
 
 $$
 C_i = \sup_m \frac{d}{|\alpha_i^{(m)}|}|| \alpha^{(m)} ||,
 $$
 
 where $\alpha^{(m)}$ is the vector whose coordinates are the $\alpha_i^{(m)}$.
\end{proof}

We are now able to conclude this part of the proof.
There are two kinds of recurrent blocks of size $n$.
Those that do not intersect any hyperplane in $\Pi$ and those that do.
Let us call $R^1_X (n)$ the number of blocks of the first kind and $R^2_X
(n)$
the number of blocks of the second kind. Lemma~\ref{lem:nonintersecting} shows that $R^1_X (n)$ is less than the constant 
$2^t\Bbeta^d$. To bound $R^2_X (n)$, we make our reasoning only
on the recurrent blocks $X_{\mathbf{y}_i,n}$ such that 
$\mathbf{y}_i$ belongs to a given domain $D$ and $\mathcal{C}(\mathbf{y}_i,n)$ intersects each of
the hyperplanes $\pi_1,\ldots,\pi_j$ and no other  (recall that there is only a
constant number of such cases to take into account). Thanks to Lemma~\ref{lem:11} and Remark~\ref{rem:11},  
moving such a block in any direction parallel to  $W$ leads to a number of pairwise distinct blocks bounded by a constant. On the other hand, when moving in a direction spanned by $\mathbf{e}_1,\ldots,\mathbf{e}_{d-\ell}$, Lemma~\ref{lem:12} shows that we can have up to $\mathcal{O}(n^{d-\ell})$ pairwise distinct blocks of size $n$. Recall from Lemma~\ref{lem:dimW} that $\ell=\dim W$ belongs to $\{1,\ldots,d-1\}$.

Furthermore, as suggested by the following example, the bound on
$R_X^2$ is tight. Consider the hyperplanes having respectively
equation $x_1=0$, $x_2=0$, \ldots, $x_{d-\ell}=0$ defining
$2^{d-\ell}$ regions of the kind $$\{(x_1,\ldots,x_d)\mid \forall
i\in\{1,\ldots,d-\ell\}: x_i \square_i 0 \}$$ with
$\square_1,\ldots,\square_{d-\ell}\in\{\ge,<\}$. Let us define a set
$X\subset\ZZ^d$ such that $\mathbf{x}=(x_1,\ldots,x_d)\in X$ if and
only if $x_1\cdots x_{d-\ell}\ge0$, i.e., if the components of
$\mathbf{x}$ are non-zero, the fact that $\mathbf{x}$ belongs to $X$
depends inly on the parity of the number of negative components amonst
the $d-\ell$ first. Observe that all $\mathcal{C}(\mathbf{v},n)$
where $\mathbf{v}=(v_1,\ldots,v_n)$ are such that, for all
$i\in\{1,\ldots,d-\ell\}$, $-n<v_i\le 0$ give different blocks.
Notice also that there is no constraint on $v_{d-\ell+1},\ldots,v_n$,
modifying any of these components gives the same block and thus such
blocks occur infinitely often. Consequently, we get a number of
distinct recurrent block equal to $n^{d-\ell}$.

\end{proof}

\section{Proof of the converse}
\label{sec:proofmain2}

In what follows, as explained in Remark~\ref{rema:word-set}, we make
no distinction between a subset of $\ZZ^d$ and its characteristic word
which belongs to $\{0,1\}^{\ZZ^d}$. We consider indifferently
complexity functions of both subsets of $\ZZ^d$ and binary words indexed by
$\ZZ^d$.

\subsection{An interesting lemma}

The following lemma based on the pigeonhole principle is of particular importance in what follows.
It is a first step to provide us with a set of relatively small local periods.
The next step is Corollary \ref{coro:lemmeQZ} which is more explicit respectively to our problem.
This lemma is proven in \cite{Quas&Zamboni:2004} for $d=2$, we simply adapt
it to our setting.

\begin{lemma}\label{lemmeQZ}
    Let $M\subset \ZZ^d$ and $C>0$ be such that $R_M (n) \leq Cn^{d-1}$ for all $n$. 
    Let $n\in \mathbb{N}$. 
    There exists $m_0\in \NN$ such
    that, for all $m$ satisfying $m<n$ and $m^d - Cn^{d-1}\geq 1$, and all
    $\z=(\z_i)_{1\leq i\leq d}$, with $||\z||\ge m_0+m$, there exists a non-zero vector
    $\v=(\v_i)_{1\leq i\leq d}$ with $||\v||\le m$, verifying
    $$M_{\z-\mathbf{v},n-m}=M_{\z,n-m} = M_{\z+\mathbf{v},n-m}.$$
\end{lemma}

\begin{proof}
    Far enough from the origin, every block is recurrent. Formally,
    for all $n\in \mathbb{N}$  
    there exists $m_0$ such that if $||\x ||\geq
    m_0$, then $M_{\x , n}$ is recurrent in $M$.

    Let $m$ and $\z=(\z_i)_{1\leq
      i\leq d}$ be such that $m<n$, $m^d - Cn^{d-1}\geq 1$ and $||\z||\ge m_0+m$.

    If $\mathbf{y}\in [\![-m+1,0]\!]^d$, then $||\z+\mathbf{y}||\ge
    m_0$. Therefore, for all $\mathbf{y}\in [\![-m+1,0]\!]^d$, the
    blocks $M_{\z+\mathbf{y},n}$ are recurrent.
    Since $R_M(n)<m^d$, by the pigeonhole principle, there exist two
    distinct integer vectors
    $\mathbf{y},\mathbf{y}'\in[\![-m+1,0]\!]^d $ such that
    $M_{\z+\mathbf{y},n} = M_{\z+\mathbf{y}',n}$.  
    Let $\mathbf{v}=\mathbf{y}-\mathbf{y}'$.  
    We observe that $||\mathbf{v}||\le m$. 
    Now consider $\mathbf{x}$ in
    $\{\z\}+[\![0,n-m-1]\!]^d$, then $\mathbf{x}$ and
    $\mathbf{x}+\mathbf{v}$ belong to
    $\{\z+\mathbf{y}\}+[\![0,n-1]\!]^d$, and, $\mathbf{x}$ and
    $\mathbf{x}-\mathbf{v}$ belong to
    $\{\z+\mathbf{y}'\}+[\![0,n-1]\!]^d$. The situation is depicted in
    Figure~\ref{fig:14}.
    \begin{figure}[htbp]
        \centering
        \includegraphics{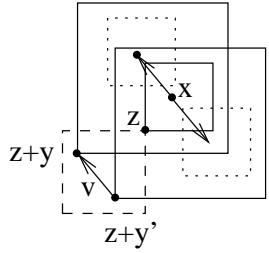}
        \caption{Sketch of the situation in Lemma~\ref{lemmeQZ}.}
        \label{fig:14}
    \end{figure}
    The conclusion follows from the fact that $M_{\z+\mathbf{y},n} =
    M_{\z+\mathbf{y}-\mathbf{v},n}$ and from the overlap of the two
    blocks.
\end{proof}

The following corollary shows that Lemma~\ref{lemmeQZ} is ``almost'' sufficient to complete the proof of Theorem~\ref{theo:main}.

\begin{coro}
\label{coro:lemmeQZ}
Let $M\subseteq\ZZ^d$.
Suppose there exists $C>0$ such that $R_M
    (n) \leq Cn^{d-1}$ for all $n$. 
Then, for all large enough $K$ there exists $L\geq 0$ such that

$$
\forall \mathbf{x}\in\ZZ^d, ||\mathbf{x}||\ge L,
    \exists \mathbf{v}, ||\mathbf{v}||\leq (2C)^{\frac{1}{d}}(5K)^{\frac{d-1}{d}}:  M \text{ is }\mathbf{v}\text{-periodic
      inside }\B(\mathbf{x},K) .
$$

\end{coro}

\begin{proof}
Suppose $d=1$: $R_M (n) \leq C$ for all $n$.
Let $x\in \{ 0,1 \}^\mathbb{Z}$ be the characteristic word of $M$.
It suffices to show that there exists three words $u$, $v$ and $w$ such that $x=\cdots uuu wvvv\cdots$.
We can of course suppose $C$ is an integer.
There exists a positive integer $n_0$ such that all words of length $C+1$ appearing in $x^+=x_{n_0}x_{n_0+1}\cdots$ and $x^- = \cdots x_{-n_0-1}x_{-n_0}$ are recurrent (i.e., appear infinitely many times in these sequences).
Thus in $x^+$ and $x^-$ the number of words of length $C$ is less than $C$.
Thus, due to Morse--Hedlund theorem (Theorem \ref{th:morse-hedlund}) both sequences are ultimately periodic: $x^+ = w^+ vvv\cdots$ and $x^- = \cdots uuuw^-$ where $u$ and $v$ are non-empty words.
Moreover it is classical to deduce from the proof of Morse--Hedlund theorem that $|u|$ and $|v|$ are both less than $C$.
This achieves the proof.

Thus we suppose $d\geq 2$.
Consider the following maps $\alpha:n\mapsto (1+Cn^{d-1})^{1/d}$, $\gamma:n\mapsto (2C)^{1/d}n^{\frac{d-1}{d}}$ and $\beta:n\mapsto (n-\gamma (n))/4$.
Let $n_0$ be such that 

\begin{itemize}
\item
$\beta$ is increasing on $[n_0 ,+\infty [$;
\item
$\gamma (n) <n $ for all $n\geq n_0$;
\item
$\gamma (n) \leq \beta (n-1)$ for all $n\geq n_0$;
\item
$n\leq 5\beta (n-1)$ for all $n\geq n_0$.
\end{itemize}

Let $K,n$ be integers such that $K\geq \beta (n_0)$, $n\geq n_0+1$ and $\beta (n-1) \leq K \leq \beta (n)$.
Let $m_0$ be given by Lemma~\ref{lemmeQZ}.
Let $m$ be an integer satisfying $\alpha (n) \leq m \leq \gamma (n) $.
Notice that $m$ satisfies the assumption of Lemma~\ref{lemmeQZ}.
We set $L = m_0+m+2K$. 
Let $\x$ satisfying $||\mathbf{x}||\geq L$.
Let $\1$ be the vector consisting of ones.
From Lemma~\ref{lemmeQZ}, since $||\mathbf{x}-2K\1||\ge m_0+m$, there exists a vector $\mathbf{v}$, with $|| \mathbf{v}||\leq m$ such that 
$M_{\x-2K\1,n-m} = M_{\x-2K\1+\mathbf{v},n-m}$.
But as $n-m\geq n-\gamma (n) = 4\beta (n) \geq 4K$,
we also have 

$$
M_{\x-2K\1,4K} = M_{\x-2K\1+\mathbf{v},4K}.
$$

Thus, $M$ is $\v$-periodic inside $\B(\mathbf{x},K)$.
Moreover, from the choice of $n_0$ we obtain that 

$$
||\v||\leq m \leq (2C)^{1/d}n^{\frac{d-1}{d}} \leq (2C)^{1/d} (5\beta (n-1))^{\frac{d-1}{d}} \leq (2C)^{1/d}(5K)^{\frac{d-1}{d}} .
$$
This completes the proof. 
\end{proof}

In order to obtain local periodicity as given in Definition~\ref{def:localp} leading to a  proof of Theorem~\ref{theo:main}, it would be nice that Corollary~\ref{coro:lemmeQZ} implies the existence of a set $V\subset \mathbb{Z}^d$  such that  $K>\sum_{\mathbf{v}\in V} ||\mathbf{v}||$.
This is the case in particular for $d=2$ (as presented below).
For $d\geq 3$ we are not able to give such a direct proof as for $d=2$.
Nevertheless this corollary is interesting because it provides us with a finite set of quite small local periods. 
We will use this when $d\geq 3$.

\subsection{Proof of Theorem~\ref{theo:main} for $d=2$}
The following two lemmata are also true in higher dimensions.
The first one is another way to settle down the underlying idea of Lemma~\ref{lemmeQZ}.
Let $\1$ be the vector consisting of ones.

\begin{lemma}
\label{lemma:distinctblocks}
Let $M\subset \mathbb{Z}^2$.
Suppose $M$ is $\mathbf{v}$-periodic inside $\mathcal{B} (\mathbf{x} , n)$, for some $\mathbf{v}\not = 0$, $||\mathbf{v}||<n$, and is not $\mathbf{w}$-periodic for any $\mathbf{w}$ with $||\mathbf{w}||<||\mathbf{v}||$.
Then,  the $||\mathbf{v}||^2$ blocks $M_{\mathbf{x} -n\1 - \mathbf{z} , 2n+||\mathbf{v}||}$ with $\mathbf{z}\in [\![ 0, ||\mathbf{v}||-1 ]\!]^2$ are pairwise distinct.
\end{lemma}

\begin{proof}
Suppose there exist two distinct vectors $\mathbf{z}_1$ and $\mathbf{z}_2$ in $[\![ 0, ||\mathbf{v}||-1 ]\!]^2$ such that $M_{\x -n\1 - \mathbf{z}_1 , 2n+||\mathbf{v}||}$ and $M_{\x -n\1 - \mathbf{z}_2 , 2n+||\mathbf{v}||}$ are equal. 
Let $\mathbf{w} = \z_1 -\z_2$. 
Then, $||\mathbf{w}||<||\mathbf{v}||$ and $M$ is $\mathbf{w}$-periodic inside $\mathcal{B} (\mathbf{x} , n)$.
\end{proof}

\begin{lemma}
\label{lemme:vordonne}
Let $n\in \NN$, $M\subset \mathbb{Z}^2$ and $\mathbf{v}_1, \dots , \mathbf{v}_k, \x_1 \dots , \x_k$ be vectors of $\mathbb{Z}^2$ such that 

\begin{enumerate}
\item 
$\mathbf{0}, \mathbf{v}_1, \dots , \mathbf{v}_k$ are pairwise distinct;
\item
$ \ell=\max  ||\mathbf{v}_i|| \leq n$;
\item
for all $i$, $M$ is $\mathbf{v}_i$-periodic inside $\mathcal{B} (\x_i , n)$ and is not $\mathbf{w}$-periodic inside $\mathcal{B} (\x_i , n)$ for any $\mathbf{w}$ with $||\mathbf{w}||<||\mathbf{v}_i||$;
\item
\label{lemme:enum:4}
for all $i$, $M$ is not $\mathbf{v}_j$-periodic inside $\mathcal{B} (\x_i , n)$ for $j<i$.
\end{enumerate}

Then, we have

$$
\# \left\{ M_{\x_i-n\1 -\mathbf{z} , 2n+\ell} \mid \z\in [\![ 0, ||\mathbf{v}_i||-1]\!]^2  , 1\leq i \leq k \right\} 
=
\sum_{i=1}^k ||\mathbf{v}_i||^2 .
$$
\end{lemma}

\begin{proof}
From Lemma~\ref{lemma:distinctblocks}, for all $i$, we have 

$$
\# \left\{ M_{\x_i - n\1 -\mathbf{z} , 2n+||\mathbf{v}_i||} \mid \mathbf{z}\in [\![ 0, ||\mathbf{v}_i||-1]\!]^2 \right\} 
= 
||\mathbf{v}_i||^2 
$$

and thus, when considering possibly larger blocks, we get 

\begin{align*}
\# \left\{ M_{\mathbf{x}_i - n\1 -\mathbf{z} , 2n+\ell} \mid \mathbf{z} \in [\![ 0, ||\mathbf{v}_i||-1]\!]^2 \right\}
= 
||\mathbf{v}_i||^2 .
\end{align*}

Let $i,j$ be such that $1\leq j < i\leq k$. It is sufficient to prove that, for all vectors $\mathbf{z}_i \in [\![ 0,||\mathbf{v}_i||-1]\!]^2$ and 
$\mathbf{z}_j \in [\![ 0, ||\mathbf{v}_j||-1]\!]^2$, 
the blocks $M_{\mathbf{x}_i -n\1 - \mathbf{z}_i , 2n+\ell}$ and $M_{\mathbf{x}_j -n\1 - \mathbf{z}_j , 2n+\ell}$ 
are distinct.

Indeed, suppose there exist $\mathbf{z}_i $ and $\mathbf{z}_j $ such that $M_{\mathbf{x}_i -n\1 - \mathbf{z}_i , 2n+\ell}$ and $M_{\mathbf{x}_j -n\1 - \mathbf{z}_j , 2n+\ell}$ are equal.
Then, since $\mathcal{B} (\x_i , n )$ is included in $\mathcal{B} (\mathbf{x}_i -n\1 - \mathbf{z}_i , 2n+\ell )$, $M$ would be $\mathbf{v}_j$-periodic inside $\mathcal{B} (\x_i , n )$.
This contradicts our assumption \eqref{lemme:enum:4}.
\end{proof}

Let us conclude with the proof of Theorem~\ref{theo:main} for $d=2$.
Let $M\subset \mathbb{Z}^2$ such that for some $C$ we have $R_M (n) \leq Cn$, for all $n$, and having all its sections definable in $\langle\ZZ;<,+\rangle$.
Thanks to Theorem~\ref{theo:muchnik}, it suffices to prove that $M$ is locally periodic.

Taking in Corollary~\ref{coro:lemmeQZ} a large enough $K$, $M$ satisfies the following property:

\medskip

{\bf (P)}  there exist two positive integers $K$ and $L'$ such that 

\begin{enumerate}
\item
\label{enum:vperiodic}
$\forall \mathbf{x}\in\ZZ^2, ||\mathbf{x}||\ge L',
    \exists \mathbf{v}, ||\mathbf{v}||\leq \sqrt{10CK} : M \text{ is }\mathbf{v}\text{-periodic
      inside }\B(\mathbf{x},K)$;
\item
\label{proof:enum:qz}
$\sqrt{10CK} \leq K$;
\item
$ 8\log (K)^3 + \frac{3CK}{\log (K)} < K$.
\end{enumerate}

Let $L''$ be such that if $||\x||\geq L''$, then the block $M_{\x , 3K}$ is recurrent.
Let $L$ be greater than $2K+L''+L'$.
We need the following lemma.

\begin{lemma}
\label{lemma:technique}
There exist vectors $\x_1 , \dots , \x_k\in \mathbb{Z}^2$ and $V = \{ \mathbf{v}_1, \dots , \mathbf{v}_k\} \subset \mathbb{Z}^2$ such that

\begin{enumerate}
\item
for all  $\mathbf{x}\in\ZZ^2$, $||\mathbf{x}||\geq L$,
there exists $\mathbf{v}\in V$ such that $ M $ is $\mathbf{v}$-periodic
inside $\B(\mathbf{x},K)$;
\item
for all $i$, $||\mathbf{v}_i || \leq \sqrt{10CK}$;
\item
$\mathbf{0}, \mathbf{v}_1, \dots , \mathbf{v}_k$ are pairwise distinct;
\item
$||\mathbf{x}_i|| \geq L $ for all $i$;
\item
for all $i$, $M$ is $\mathbf{v}_i$-periodic inside $\mathcal{B} ( \x_i , K)$ and is not $\mathbf{w}$-periodic inside $\mathcal{B} ( \x_i  , K)$ for any $\mathbf{w}$ with $||\mathbf{w}||<||\mathbf{v}_i||$;
\item
for all $i$, $M$ is not $\mathbf{v}_j$-periodic inside $\mathcal{B} ( \x_i , K)$ for $j<i$.
\end{enumerate}
\end{lemma}

\begin{proof}
Let $(\mathbf{u}_i)$ be a sequence consisting of all vectors of $\mathbb{Z}^2\setminus \{ \mathbf{0} \}$, appearing only once, which is non-decreasing with respect to their norms.
Let $RB$ be the set of recurrent blocks of size $K$.
Note that from Property {\bf (P)} all elements of $RB$ are $\mathbf{v}$-periodic for some $\mathbf{v}$ whose norm is less than $\sqrt{10CK}$. 
Let $RB_0$ be the set of recurrent blocks of size $K$ having $\mathbf{u_0}$ as a local period.
Of course, $RB_0$ can be empty. 
Let $RB_1$ be the subset of $RB \setminus RB_0$ whose blocks have the local period $\mathbf{u_1}$.
Observe that these blocks do not have $\mathbf{u_0}$ as a local period.
Continuing this way, we obtain finitely many non-empty subsets $RB_{i_1}, \dots ,RB_{i_k}$ of $RB$ such that 

\begin{enumerate}
\item
$RB = \cup_{n=1}^k RB_{i_n}$;
\item
the blocks of $RB_{i_n}$ have $\mathbf{u}_{i_n}$ as a local period;
\item
for all $s\in \{ 1,\dots , k\}$ and $j<i_s$, $\mathbf{u}_{j}$ is not a local period for blocks in $RB_{i_s}$.
\end{enumerate}

Property {\bf (P)} ensures that for all $s$, $||\mathbf{u}_{i_s}||\leq \sqrt{10CK}$.
We set $\mathbf{v}_s = \mathbf{u}_{i_s}$.
For all $s$, there exists $\mathbf{x}_s \in \mathbb{Z}^2$ such that $M$ is $\mathbf{v}_s$-periodic inside $\mathcal{B} (\mathbf{x}_s , K)$.
The choice of $L$ allows us to suppose $||\mathbf{x}_s || \geq L$. 
This concludes the proof.
\end{proof}

Let $\ell = \max_{\mathbf{v}\in V} ||\mathbf{v}||$. 
Lemma~\ref{lemma:technique} provides us with vectors $\x_1 , \dots , \x_k\in \mathbb{Z}^2$ and a set $V = \{ \mathbf{v}_1, \dots , \mathbf{v}_k\} \subset \mathbb{Z}^2$ fulfilling the hypothesis of Lemma~\ref{lemme:vordonne}.
To get local periodicity it remains to show that $\sum_{\mathbf{v}\in V} ||\mathbf{v}|| < K$.
Applying Lemma~\ref{lemme:vordonne}, we obtain

$$
\sum_{\mathbf{v}\in V} ||\mathbf{v}||^2 \leq R_M (2K+\ell) \leq C(2K+\ell) .
$$

Using Property {\bf (P)}, we deduce that $\ell\leq \sqrt{10CK} \leq K$. 
Hence 

$$
\sum_{\mathbf{v}\in V} ||\mathbf{v}||^2 \leq 3CK  .
$$

Let $\delta_n = \# \{ \v\in V \mid ||\mathbf{v}||=n \}$.
Then, we have

$$
 \log (K) ^2 \sum_{n=\lfloor\log (K)\rfloor}^{\lfloor\sqrt{10CK}\rfloor} \delta_n
 \leq
 \sum_{n=1}^{\lfloor\sqrt{10CK}\rfloor} \delta_n n^2 
 =
 \sum_{\mathbf{v}\in V} ||\mathbf{v}||^2
 \leq 3CK .
$$

Consequently, we get

$$
\sum_{n=\lfloor\log (K)\rfloor}^{\lfloor\sqrt{10CK}\rfloor} \delta_n \leq \frac{3CK}{\log (K)^2} . 
$$

Hence, using Cauchy--Schwarz inequality, we deduce that

\begin{align*}
\sum_{n=\lfloor\log (K)\rfloor}^{\lfloor\sqrt{10CK}\rfloor} \delta_n n 
= &
\sum_{n=\lfloor\log (K)\rfloor}^{\lfloor\sqrt{10CK}\rfloor} \sqrt{\delta_n}  \sqrt{\delta_n} n \\
\le & 
\left(
\sum_{n=\lfloor\log (K)\rfloor}^{\lfloor\sqrt{10CK}\rfloor} \delta_n 
\right)^{1/2}
\left(
\sum_{n=\lfloor\log (K)\rfloor}^{\lfloor\sqrt{10CK}\rfloor} \delta_n n^2 
\right)^{1/2}
\leq
 \frac{3CK}{\log (K)} .
\end{align*}

Using the fact that, for all $n$, $$\delta_n \leq \# \{ \v \in \mathbb{Z}^2 \mid ||\mathbf{v}||=n \}=\# \bigcup_{i=-n}^n\{(i,n), (i,-n), (n,i), (-n,i)\}=8n,$$ we deduce that

\begin{align*}
\sum_{\v \in V} ||\mathbf{v}|| 
& = \sum_{n=1}^{\lfloor\sqrt{10CK}\rfloor} \delta_n n \leq 
\sum_{n=1}^{\lfloor\log (K)\rfloor-1} \delta_n n + \frac{3CK}{\log (K)}\\
& \leq 
8\sum_{n=1}^{\lfloor\log(K)\rfloor-1}  n^2 + \frac{3CK}{\log (K)} \leq 8\log (K)^3 + \frac{3CK}{\log (K)} <K .
\end{align*}

This concludes the proof for $d=2$.

\medskip

Observe that this kind of computation gives nothing in dimension
$d\geq 3$ because, given some constants $a,b,c$ one should find some
sufficiently large $K$ and a function $\alpha (K)$ such that $a +
b\alpha (K)^{d+1} + \frac{cK^{d-1}}{\alpha (K)} <K$. This is clearly
not possible. But, maybe some of the above inequalities could be
improved in order to have a direct proof for all $d\geq 3$.  We leave
this as an open question.

In fact, as we did not used the assumption that each section of $M$ is definable, we prove more than what was expected.
We have shown that the hypothesis $R_M (n)\in \mathcal{O} (n)$ implies the local periodicity of $M$.
Thus, the definability of the sections is only needed to go from the local periodicity of $M$ to the definability of $M$.
We do not know whether it is true in higher dimensions.
This is related to the previous remark.

\subsection{Proof of Theorem~\ref{theo:main} in the general case}
It remains to prove that if $M$ is a subset in $\ZZ^d$ such that $R_M(n)\in\mathcal{O}(n^{d-1})$ and every section is definable in $\langle\ZZ;<,+\rangle$, then it is definable in $\langle\ZZ;<,+\rangle$.
Let us first recall some crucial results of Muchnik. 
For all $A\subset \mathbb{Z}^d$ and $\mathbf{v}\in \mathbb{Z}^d$, we define the {\em border of $A$ in the direction $\mathbf{v}$} as 
$$\mathbf{Bd} (A,\mathbf{v}):= \{ \mathbf{x}\in A \mid \mathbf{x}+\mathbf{v}\not \in A \}.$$
    \begin{figure}[htbp]
        \centering
        \includegraphics{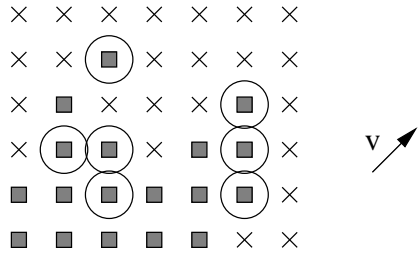}
        \caption{Illustration of the border of $A$ in the direction $\mathbf{v}$.}
        \label{fig:new03}
    \end{figure}
In Figure~\ref{fig:new03}, points in $A$ are squares and points belonging to $\mathbf{Bd} (A,\mathbf{v})$ are inside circles.

The following two lemmata show that it is sufficient to prove that, for some $\mathbf{v}$, the set $\mathbf{Bd}(M,\mathbf{v})$ is locally periodic.

\begin{lemma}
\label{lemma:definablesections}\cite[Lemma~1.0]{Muchnik:2003}
Let $A\subset \mathbb{Z}^d$ and $\mathbf{v}\in \mathbb{Z}^d$.
If all sections of $A$ are definable, then all sections of 
$\mathbf{Bd}(A,\mathbf{v})$ are definable. 
\end{lemma}

\begin{lemma}\cite[Lemma~1.2]{Muchnik:2003}
\label{lemma:definablewrtbd}
Let $A\subset \mathbb{Z}^d$ and $\mathbf{v}\in \mathbb{Z}^d$.
The set $A$ is definable in terms of $\mathbf{Bd}(A,\mathbf{v})$, $\mathbf{Bd}(A,-\mathbf{v})$ and a finite number of sections of $A$, i.e., 
$A$ can be defined by a formula involving addition, order and unary predicates for the sets $\mathbf{Bd}(A,\mathbf{v})$, $\mathbf{Bd}(A,-\mathbf{v})$ and for a finite number of sections of $A$.
\end{lemma}

The next lemma suggests to proceed by induction on the cardinality of the set of local periods to prove that $\mathbf{Bd}(M,\mathbf{v})$ is locally periodic.

\begin{lemma}\cite[Lemma~1.1]{Muchnik:2003}
\label{lemma:V-vperiodic}
Let $A\subset \mathbb{Z}^d$.
Suppose there exist $K$, $L\geq 0$ and a finite set $V\subset \mathbb{Z}^d\setminus \{ \mathbf{0}\}$ verifying
that, for all $\mathbf{x}\in\ZZ^d$ with $||\mathbf{x}||\ge L$,
the set $A$ is $V$-periodic inside
$\B(\mathbf{x},K)$. 
Then, for all $\mathbf{v}\in \mathbb{Z}^d$ and all $\mathbf{x}\in\ZZ^d$ with $||\mathbf{x}||\ge L$,
the set $\mathbf{Bd}(A,\mathbf{v})$ is $V\setminus \{ \mathbf{v} , -\mathbf{v}\}$-periodic inside $\B(\mathbf{x},K)$.
\end{lemma}

To complete the proof of Theorem~\ref{theo:main} we will proceed by induction on the cardinality of $V$. 

We will say that a set $M\subseteq\ZZ^d$ satisfies ${\rm Hyp} (k)$ whenever :

\begin{enumerate}
\item
there exists $C>0$ such that, for all large enough $n\in\NN$, $R_{M}
    (n) \leq Cn^{d-1}$; 
\item
all sections of $M$ are definable; 
\item
\label{Pk:cond3}
there exist $K$, $L\geq 0$ and $V\subset \mathbb{Z}^d\setminus \{ \mathbf{0}\}$ verifying
\begin{enumerate}
\item
$\forall \mathbf{x}\in\ZZ^d, ||\mathbf{x}||\ge L,
    \exists \mathbf{v}\in V: M \text{ is }\mathbf{v}\text{-periodic
      inside }\B(\mathbf{x},K)$,
\item
$\forall \mathbf{v}\in V$, $||\mathbf{v}||\leq (2C)^{\frac{1}{d}}(5K)^{\frac{d-1}{d}}< K$, and,
\item
$\# V = k$.
\end{enumerate} 
\end{enumerate}

We note that, from Corollary~\ref{coro:lemmeQZ}, 
the subsets of $\ZZ^d$ such that $R_M(n)\in\mathcal{O}(n^{d-1})$ and having all its sections definable
        in $\langle\ZZ;<,+\rangle$, satisfy ${\rm Hyp} (k)$ for some $k$.
Hence, to complete the proof of Theorem~\ref{theo:main} it is sufficient to prove that the following assertion $P(k)$ is true for all $k$.
$$P(k):\quad \text{ If }M\subseteq\ZZ^d \text{ satisfies Hyp(}k\text{), then it is definable}.$$

If $M$ satisfies ${\rm Hyp} (1)$ then \eqref{Pk:cond3} corresponds to the local periodicity of $M$.
Hence $P(1)$ is true. 

Let $k\ge 1$ be such that $P(1), \dots , P(k)$ are true. 
Let $M$ be a subset of $\ZZ^d$ satisfying ${\rm Hyp} (k+1)$.
It remains to prove that $M$ is definable.
Let $C$, $V$, $K$ and $L$ be given by ${\rm Hyp} (k+1)$ for $M$.

\begin{lemma}
\label{lemma:epsilon-comp}
For all $\mathbf{v}\in \mathbb{Z}^d$ and $\epsilon >0$ there exists $n_0$ such that if $n\geq n_0$, then $R_{\mathbf{Bd}(M,\mathbf{v})  } (n) \leq (C+\epsilon ) n^{d-1}$.
\end{lemma}

\begin{proof}
It suffices to see that $R_{\mathbf{Bd}(M,\mathbf{v})  } (n) \leq R_{M } (n+||\mathbf{v}||)$ for all $n\in \mathbb{N}$.
\end{proof}

\begin{lemma}
\label{lemma:hypk}
For all $\mathbf{v}\in V$ the sets $\mathbf{Bd}(M,\mathbf{v})$ and $\mathbf{Bd}(M,-\mathbf{v})$ satisfy ${\rm Hyp} (k)$.
\end{lemma}

\begin{proof}
Let $\mathbf{w}\in \{ \mathbf{v}, -\mathbf{v}\}$.
From Lemma~\ref{lemma:epsilon-comp} there exists $\epsilon >0$ such that 

$$
(2(C+\epsilon ))^{\frac{1}{d}}(5K)^{\frac{d-1}{d}}< K \hbox{ and } R_{\mathbf{Bd}(M,\mathbf{w})  } (n) \leq (C+\epsilon ) n^{d-1}
$$ 

for all large enough $n$. 
From Lemma~\ref{lemma:definablesections} all sections of $\mathbf{Bd}(M,\mathbf{w})$ are definable.
We take for $\mathbf{Bd}(M,\mathbf{w})$ the constants $K,L$ and the set $V$ of $M$.
With Lemma~\ref{lemma:V-vperiodic} we know that $\mathbf{Bd}(M,\mathbf{w})$ is $V\setminus \{ \mathbf{v} , -\mathbf{v}\}$-periodic inside $\B(\mathbf{x},K)$ for all $\mathbf{x}$ with $||\mathbf{x}||\geq L$, and, thus satisfies ${\rm Hyp} (k)$.
\end{proof}

We can now conclude the proof of Theorem~\ref{theo:main}.
Let $\mathbf{v}\in V$.
From Lemma~\ref{lemma:hypk}, $\mathbf{Bd}(M,\mathbf{v})$ and $\mathbf{Bd}(M,-\mathbf{v})$ satisfy ${\rm Hyp} (k)$.
BY induction hypothesis both $\mathbf{Bd}(M,\mathbf{v})$ and $\mathbf{Bd}(M,-\mathbf{v})$ are definable. 
Hence, by Lemma~\ref{lemma:definablewrtbd}, $M$ is definable.

\medskip

We recall that for $d=2$ we do not need the definability of the sections to deduce the local periodicity. 
Observe that for dimensions $d\geq 3$ it is not the case: We need the definability of the sections to show that the boarders are definable and then, using the induction process and the Muchnik criterion, to conclude the proof. 

\section{Related works}
\label{section:comments}

In this section we recall some well-known works on the relations between the block complexity and the periodicity for subsets of $\ZZ^d$ and $\RR^d$.
We end this section with some comments, related to these works, on our main result with the additional hypothesis of repetitiveness.

\subsection{Delone sets and repetitiveness (or recurrence)}

Let us recall some terminology in \cite{Lagarias&Pleasants:2002}.
Let $X$ be a subset of  $\RR^d$.
We say that $X$ is a $(r,R)$-{\em Delone set} if it has the following two properties:

\begin{enumerate}
\item
{\em Uniform Discreteness}. 
Each open ball of radius $r$ in $\RR^d$ contains at most one point of $X$.
\item
{\em Relative density}.
Each closed ball of radius $R$ in $\RR^d$ contains at least one point of $X$.
\end{enumerate}

Observe that the subsets of $\ZZ^d$ are uniformly discrete for $r=1/2$.

The {\em period lattice} of a Delone set $X$ is the lattice of translation symmetries given by $\Lambda_X = \{ \mathbf{p}\in \RR^d | X+\mathbf{p} = X \}$.
It is a free abelian group with rank between $0$ and $d$.
When the rank of $\Lambda_X$ is equal to $d$, we say that $X$ is an {\em ideal crystal}.

Let $B (\mathbf{x} , t )$ stands for the open ball centered in $\mathbf{x}\in \RR^d$ of radius $t$.
The set $\P_X (\x , t) = X\cap B (\mathbf{x} , t )$ is called a $t$-patch of $X$. 
In the sequel we will consider that two $t$-patches are equal when they are equal up to translation.

The Delone set $X$ is {\em repetitive} if for each $t$ there is $M(t)>0$ such that every closed ball $B$ of radius $t$ contains all $t$-patches of $X$ (up to translation).
It is {\em linearly repetitive} if there exists a constant $C$ such that $M(t)\leq C t$ for all $t>0$.
We denote by $p_X (t)$ the number (possibly infinite) of different $t$-patches (up to translation).

\subsection{The Period Conjecture}
\label{subsec:periodconjecture}

The {\em Period Conjecture} is stated in \cite{Lagarias&Pleasants:2002} for Delone sets in $\RR^d$. 

\medskip

{\bf Period Conjecture.}
{\em For each integer $j=1,\dots , d$, there is a positive constant $c_j (r,R)$ such that any $(r,R)$-Delone set $X$ of $\RR^d$ satisfying 

$$
p_X ( t) < c_j (r,R) t^{d-j+1} \hbox{ for all } t> t_0 (X),
$$

for some $t_0 (X)$, has $j$ linearly independent periods.
}

\medskip

Lagarias and Pleasants showed in \cite{Lagarias&Pleasants:2002} that this conjecture is true for $j=n$.

\begin{theo}
If a $(r,R)$-Delone set $X\subset \RR^d$ has a single value $t$ such that 
$$
p_X (t) < \frac{t}{2R}
$$

then $X$ is an ideal crystal.
\end{theo}

In \cite{Epifanio&Koskas&Mignosi:2003} and \cite{Quas&Zamboni:2004} it is also proven for $j=1$ and subsets of $\ZZ^2$ that are not necessarily Delone sets (see Theorem \ref{theo:QZ} below).

In \cite{Lenz:2004}, D. Lenz answers positively to a conjecture in \cite{Lagarias&Pleasants:2003} saying that with the extra hypothesis of linear repetitiveness the Period Conjecture is not far to be true for $j=1$.

\begin{theo}
Every aperiodic linearly repetitive Delone set $X\subset \RR^d$ satisfies

$$
\liminf_{t\to\infty} \frac{p_X (t)}{t^d} >0.
$$
\end{theo}

\subsection{Nivat's Conjecture}

Let us recall the {\em Nivat's Conjecture} \cite{Nivat:1997} stated in 1997 (see also \cite{Berthe&Vuillon:2000}).

\begin{conj}[M. Nivat]
    Let $M$ be a subset of $\ZZ^2$. If there exist $n_1,n_2>0$ such
    that the function $p_M$ counting the number of distinct
    $(n_1\times n_2)$-blocks occurring in $M$ is such that
    $p_M(n_1,n_2)\le n_1n_2$, then $M$ is periodic.
\end{conj}

Nivat's Conjecture cannot be an equivalence because the
converse does not hold: there exists a periodic set $M$ in $\ZZ^2$
such that $p_M(n_1,n_2)> n_1n_2$, for all $n_1,n_2$, see \cite[p.~49]{Berthe&Vuillon:2000}. 
Moreover, it cannot be true for dimensions $d$ strictly greater than $2$ as shown in \cite{Sander&Tijdeman:2000}. 
A weaker form of Nivat's Conjecture is the following. Let $\alpha<1$. If there exist $n_1,n_2>0$ such
     such that $p_M(n_1,n_2)\le \alpha n_1n_2$,
    then $M$ is periodic.
It was proven for $\alpha=1/144$ in \cite{Epifanio&Koskas&Mignosi:2003} and improved for $\alpha=1/16$ in \cite{Quas&Zamboni:2004}. This latter result is a consequence of the following one.

\begin{theo}\cite{Quas&Zamboni:2004}
\label{theo:QZ}
    Let $M$ be a subset of $\ZZ^2$. If there exist $n_1,n_2>0$ such
    that for all $(2n_1\times 2n_2)$-blocks $B$ the function counting the maximum number of distinct
    $(n_1\times n_2)$-blocks occurring in $B$ is
    less or equal to $\frac{n_1n_2}{16}$, then $M$ is periodic.
\end{theo}

As explained in Section \ref{subsec:periodconjecture}, this theorem is related to the Period Conjecture.

For a survey on the relationships existing between periodicity and block complexity in $\ZZ^d$ we refer to \cite{Cassaigne:2000} and \cite{Tijdeman:2006}. 
We can also mention \cite{Sander&Tijdeman:2002}, see below.

\subsection{The Periodicity Principle}

Let us recall some terminology introduced in
\cite{Sander&Tijdeman:2000}.  Let $S \in \{ 0,1 \}^{\ZZ^d}$.  For all
finite subset $B$ of $\ZZ^d$ we define the $B$-patterns of $S$ to be
the functions $S^{(\mathbf{v})} : B \to \{ 0,1 \}$ defined by
$S^{(\mathbf{v})} (\mathbf{x}) = S(\mathbf{x}+\mathbf{v})$,
$\mathbf{v}\in \ZZ^d$.  The set of all $B$-patterns is $P_S (B) = \{
S^{(\mathbf{v})} | \mathbf{v}\in \ZZ^d \}$.  When $B$ is a cube of
size $n$, then $P_S (B)$ is equal to the block complexity $p_S (n)$.
In \cite{Sander&Tijdeman:2000} the following conjecture is stated .

\medskip

{\bf Periodicity Principle.}
{\em 
Let $S \in \{ 0,1 \}^{\ZZ^d}$ and $B$ be a finite subset of  $\ZZ^d$.
If $\# P_S (B) \leq \# B$, then $S$ is periodic. 
}

\medskip

It is shown in \cite{Sander&Tijdeman:2000} that the Periodicity
Principle is true in dimension 1 and turns to be false in higher
dimensions without some additional assumptions.  If the Nivat's
Conjecture holds true, it could be considered as a {\em restricted}
Periodicity Principle for rectangle blocks.

\subsection{Final comments}

First, one can show that for subsets $M$ of $\ZZ^d$ having all its
blocks of size $n$ occurring infinitely many times, definability
implies periodicity.

Now let us give a corollary of our main result
(Theorem~\ref{theo:main}) in the context of repetitive subsets of
$\ZZ^d$, and then comment it relatively to the previously described
works.

\begin{theo}
\label{theo:coromain}
Let $M$ be a repetitive subset of $\ZZ^d$.
The following statements are equivalent:

\begin{enumerate}
\item
\label{hyp1}
$R_M(n)\in\mathcal{O}(n^{d-1})$ and every section is definable in $\langle\ZZ;<,+\rangle$;
\item
\label{hyp2}
for all $k \in \{ 1, \dots , d-1 \}$ every $(d-k)$-dimensional section of $M$ has a recurrent block
    complexity in $\mathcal{O}(n^{d-k-1})$;
\item 
$M$ is an ideal crystal.
\end{enumerate}
\end{theo}

\begin{proof}
It suffices to use Theorem \ref{theo:muchnik} and to see that when a set is semi-linear and repetitive then it is a ideal crystal (the converse being also true). 
\end{proof}

Observe that in the above result, assumptions \eqref{hyp1} or
\eqref{hyp2} are much stronger than those of the Period Conjecture.
Nevertheless we obtain a necessary and sufficient condition to be an
ideal crystal in terms of block complexity.  Moreover, for $j=2$, the
Period Conjecture expects only two linearly independent periods when
we obtain an ideal crystal, that is $d$ linearly independent periods.
 
Similar comments are also valid for the Periodicity Principle.  Thus,
we obtain a restricted Periodicity Principle as it holds for cubic
blocks $B$.  But we obtain a much stronger result than just a one
dimensional periodicity conclusion.

A natural question is whether Theorem~\ref{theo:coromain} remains true
for any finite subset $B$ in order to have a Periodicity Principle for
repetitive sets.  Another question is to ask whether the Periodicity
Principle or the Period Conjecture are true for some well-known
families of Delone sets or tilings like those of finite type, linearly
repetitive, substitutive, ...

\subsection*{Acknowledgements}
The authors thank the referees for their valuable and relevant comments.



\end{document}